 \documentclass{amsart}
 \usepackage{euscript}
 \usepackage[top=1in,bottom=1in,left=1.25in,right=1.25in]{geometry} 

 \usepackage{amsfonts}
 \usepackage{amsmath}
 \usepackage{amssymb}
 \usepackage{latexsym,color}
 \usepackage{amsbsy}
 \usepackage{graphicx}
 \usepackage{calc}
 \usepackage{ifthen}
 \usepackage{epsfig}
 \usepackage{pstricks,pst-node}

 \newrgbcolor{skyblue}{.53 .81 .98} 
 \newrgbcolor{firebrick1}{1 .19 .19}

 \theoremstyle{plain} 
 \newtheorem{thm}{Theorem}[section]
 \newtheorem{con}{Conjecture}[section]
 \newtheorem{lemma}[thm]{Lemma}
 \newtheorem{cor}[thm]{Corollary}
 \newtheorem{pro}[thm]{Proposition}

 \theoremstyle{definition} 
 \newtheorem{example}[thm]{Example}
 \newtheorem{definition}[thm]{Definition}
 \newtheorem{remark}[thm]{Remark}
 \newtheorem{Algorithm}[thm]{Algorithm}

 \newcommand{\ncom}{\newcommand}
 \ncom{\la}{\lambda}
 \ncom{\bm}{\boldmath}
 \ncom{\noi}{\noindent}
 \ncom{\bq}{\begin{equation}}
 \ncom{\eq}{\end{equation}}
 \ncom{\beqn}{\begin{eqnarray*}}
 \ncom{\eeqn}{\end{eqnarray*}}
 \ncom{\ba}{\begin{array}}
 \ncom{\ea}{\end{array}}
 \ncom{\beq}{\begin{eqnarray}}
 \ncom{\eeq}{\end{eqnarray}}
 \ncom{\nno}{\nonumber}
 \ncom{\hs}{\mbox{\hspace{.25cm}}}
 \ncom{\rar}{\rightarrow}
 \ncom{\Rar}{\Rightarrow}
 \ncom{\noin}{\noindent}
 \ncom{\bc}{\begin{center}}
 \ncom{\ec}{\end{center}}
 \ncom{\sz}{\scriptsize}
 \ncom{\fpd}{\Phi(\pi^{'})}
 \ncom{\fp}{\Phi(\pi) }
 \ncom{\nk}{\left< \begin{array}{c}
                        n\\k \end{array} \right>}
 \ncom{\nd}{1^{'},2^{'},\cdots,n^{'}}
 \ncom{\R}{I\!\!R}
 \ncom{\de}{\bigtriangleup (F_{2n}, \leq)}
 \ncom{\del}{\bigtriangleup}
 \ncom{\cov}{<\!\!\!\!\cdot }
 \ncom{\bt}{\begin{thm}}
 \ncom{\bcon}{\begin{con}}
 \ncom{\et}{\end{thm}}
 \ncom{\econ}{\end{con}}
 \ncom{\bl}{\begin{lemma}}
 \ncom{\el}{\end{lemma}}
 \ncom{\bco}{\begin{cor}}
 \ncom{\ds}{\displaystyle}
 \ncom{\eco}{\end{cor}}
 \ncom{\bp}{\begin{pro}}
 \ncom{\ep}{\end{pro}}
 \ncom{\bex}{\begin{example}}
 \ncom{\eex}{\end{example}}
 \ncom{\bd}{\begin{definition}}
 \ncom{\ed}{\end{definition}}
 \ncom{\brm}{\begin{remark}}
 \ncom{\erm}{\end{remark}}
 \ncom{\bal}{\begin{Algorithm}}
 \ncom{\eal}{\end{Algorithm}}

 \ncom{\pf}{\begin{proof}}
 \ncom{\epf}{\end{proof}}
 \ncom{\be}{\begin{enumerate}}
 \ncom{\ee}{\end{enumerate}}
 \ncom{\s}{\subset}
 \ncom{\cc}{\mathcal{C}}
 \ncom{\cf}{\mathcal{F}}
 \ncom{\ac}{{\mathcal{A}}(G,{\mathcal{C}})}
 \ncom{\tc}{{\mathcal{T}}(G,{\mathcal{C}})}
 \ncom{\trc}{{\mathcal{TR}}(G,{\mathcal{C}})}
 \ncom{\kc}{{\mathcal{K}}(G)}
 \ncom{\kcx}{{\mathcal{K}}(G_X')}
 \ncom{\zc}{{\mathcal{Z}}(G)}
 \ncom{\mn}{\mathbb{N}}
 \ncom{\mr}{\mathbb{R}}
 \ncom{\mq}{\mathbb{Q}}
 \ncom{\stack}[2]{\genfrac{}{}{0pt}{}{#1}{#2}} 

 \ncom{\tcr}{\bf\textcolor{firebrick1}}
 \ncom{\tcb}{\bf\textcolor{skyblue}}

 \newenvironment{entry}
   {\begin{list}{}%
      {%
        \setlength{\labelwidth}{35pt}%
        \setlength{\leftmargin}{\labelwidth+\labelsep}%
      }%
   }%
   {\end{list}}

   \newlength{\Mylen}
   \newcommand{\Lentrylabel}[1]{%
     \settowidth{\Mylen}{\emph{#1}}%
     \ifthenelse{\lengthtest{\Mylen > \labelwidth}}%
        {\parbox[b]{\labelwidth}
          {\makebox[0pt][l]{\emph{#1}}\\}}%
        {\emph{#1}}
     \hfil\relax}

 \newenvironment{Lentry}%
   {%
    \begin{entry}}%
   {\end{entry}}

\begin{document}
 \title[Alternating Cones]{Cones of closed alternating walks and trails}
 \author[Bhattacharya]{Amitava Bhattacharya}
 \address{Bhattacharya: Department of Mathematics, Statistics, and Computer
 Science\\
 University of Illinois at Chicago\\
 Chicago, Illinois 60607-7045, USA\\
 Phone: (312) 413 2163\\
 Fax: (312) 996 1491}
 \email{amitava@math.uic.edu}
 \author[Peled]{Uri N. Peled}
 \address{Peled: Department of Mathematics, Statistics, and Computer
 Science\\
 University of Illinois at Chicago\\
 Chicago, Illinois 60607-7045, USA\\
 Phone: 312 413 2156
 Fax: (312) 996 1491}
 \email{uripeled@uic.edu}
 \author[Srinivasan]{Murali K. Srinivasan*}
 \address{Srinivasan: Department of Mathematics\\
 Indian Institute of Technology, Bombay\\
 Powai, Mumbai 400076, INDIA\\
 Phone: 91-22-2576 7484\\
 Fax: 91-22-2572 3480}
 \email{mks@math.iitb.ac.in}
 \thanks{UNP and MKS would like to thank Professor Martin Golumbic for his kind
 invitation to visit the Caesarea Edmond Benjamin de Rothschild
 Foundation Institute for Interdisciplinary Applications of Computer
 Science at the University of Haifa, Israel during May--June 2003, where
 part of this work was carried out. The warm hospitality and partial
 support of this visit from CRI is gratefully acknowledged.
 \\ * Corresponding author}

 \dedicatory{Dedicated to the memory of Malka Peled}
 \date{October 22, 2005; Revised October 12, 2006; minor typo January 9, 2007. To appear in Linear Algebra and Its Applications}
 \begin{abstract}
 Consider a graph whose edges have been colored red and blue. Assign
 a nonnegative real weight to every edge so that at every vertex,
 the sum of the weights of the incident red edges equals the sum of
 the weights of the incident blue edges. The set of all such
 assignments forms a convex polyhedral cone in the edge space,
 called the \emph{alternating cone}. The integral (respectively,
 $\{0,1\}$) vectors in the alternating cone are sums of
 characteristic vectors of closed alternating walks (respectively,
 trails). We study the basic properties of the alternating cone,
 determine its dimension and extreme rays, and relate its dimension
 to the majorization order on degree sequences. We consider whether
 the alternating cone has integral vectors in a given box, and use
 residual graph techniques to reduce this problem to the one of
 searching for an alternating trail connecting two given vertices.
 The latter problem, called \emph{alternating reachability}, is
 solved in a companion paper along with related results.

 \end{abstract}
 \maketitle

 \section{Introduction and Summary}\label{sec1}

 Consider a directed graph. Assign a nonnegative real weight to
 every arc so that at every vertex, the total weight of the incoming
 arcs is equal to the total weight of the outgoing arcs. The set of
 all such assignments forms a convex polyhedral cone in the arc
 space, called the cone of circulations, and is a basic object of
 study in network flow theory. For instance, placing integral upper
 and lower bounds on every arc and asking whether there is an
 integral vector in the cone of circulations meeting these bounds
 leads to Hoffman's circulation theorem (see the book~\cite{ff}).
 Now consider an undirected analog of the situation above. Take an
 undirected graph whose edges have been colored red and blue. Assign
 a nonnegative real weight to every edge so that at every vertex,
 the total weight of the incident red edges equals the total weight
 of the incident blue edges. The set of all such assignments forms a
 convex polyhedral cone in the edge space, called the alternating
 cone. In this paper and the companion paper~\cite{BPS2}, we study
 the basic theory of the alternating cone. Here we consider its
 extreme rays, integral vectors, and dimension. We also relate it to
 threshold graphs and majorization order on degree sequences. We
 reduce the problem of finding an integral vector in the alternating
 cone whose components satisfy given upper and lower bounds to the
 problem of searching for an alternating trail connecting two given
 vertices in a 2-colored graph (recall that in the directed case,
 the corresponding problem is reduced to the problem of searching
 for a directed path from one given vertex to another in a suitable
 residual directed graph). This latter problem, called alternating
 reachability, generalizes the problem of searching for an
 augmenting path with respect to a matching in a non-bipartite graph
 and is solved in~\cite{BPS2} by generalizing the blossom forest
 algorithm of Edmonds. We now give precise definitions and an
 outline of our results.

 Let $G=(V,E)$ be an undirected graph (we allow parallel edges but
 not loops). Assume that the edges of $G$ are colored red or blue,
 the coloring being given by ${\mathcal C}: E \rar \{R,B\}$. We say
 that $(G,\mathcal{C})$ is a \emph{2-colored graph}. Consider the
 real vector space ${\mathbb R}^E$, with coordinates indexed by the
 set of edges of $G$. We write an element $x \in {\mathbb R}^E$ as
 $x= (x(e) : e \in E)$. For a subset $F \subseteq E$ and $v \in V$,
 $F(v)$ denotes the set of all edges in $F$ incident with $v$. For a
 subset $F \subseteq E$, $F_R$ (respectively, $F_B$) denotes the set
 of red (respectively, blue) edges in $F$. For an edge $e \in E$,
 the characteristic vector $\chi(e) \in {\mathbb R}^E$ is defined by
 $\chi(e)(f) = \left\{ \ba{cc} 1, & \mbox{if }f=e \\
                                          0, & \mbox{if }f \neq e
                                   \ea \right.$.
 The \emph{red degree} $r(v)$ (respectively, \emph{blue degree}
 $b(v)$) of a vertex $v \in V$ is the number of red (respectively,
 blue) edges incident with $v$.

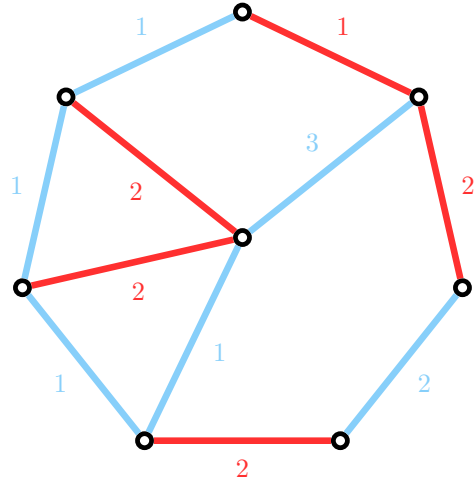
\begin{figure}
   \begin{center}
     \begin{pspicture}(-4.5,-3)(4.5,3.1)
        \psset{radius=0.125cm}
        \psset{linewidth=1.8pt}
        \SpecialCoor
        \degrees[360]
        \Cnode(0;1){n0}
        \Cnode(3;90){n1}
        \Cnode(3;!360 7 div 90 add){n2}
        \Cnode(3;!2 360 mul 7 div 90 add){n3}
        \Cnode(3;!3 360 mul 7 div 90 add){n4}
        \Cnode(3;!4 360 mul 7 div 90 add){n5}
        \Cnode(3;!5 360 mul 7 div 90 add){n6}
        \Cnode(3;!6 360 mul 7 div 90 add){n7}
        \psset{linecolor=firebrick1}
        \psset{linewidth=2.5 pt}
        \ncline{n0}{n2} \Aput{\textcolor{firebrick1}{$2$}}
        \ncline{n0}{n3} \Aput{\textcolor{firebrick1}{$2$}}
        \ncline{n4}{n5} \Bput{\textcolor{firebrick1}{$2$}}
        \ncline{n6}{n7} \Bput{\textcolor{firebrick1}{$2$}}
        \ncline{n7}{n1} \Bput{\textcolor{firebrick1}{$1$}}
        \psset{linecolor=skyblue}
        \ncline{n0}{n4} \Aput{\textcolor{skyblue}{$1$}}
        \ncline{n0}{n7} \Aput{\textcolor{skyblue}{$3$}}
        \ncline{n1}{n2} \Bput{\textcolor{skyblue}{$1$}}
        \ncline{n2}{n3} \Bput{\textcolor{skyblue}{$1$}}
        \ncline{n3}{n4} \Bput{\textcolor{skyblue}{$1$}}
        \ncline{n5}{n6} \Bput{\textcolor{skyblue}{$2$}}
    \end{pspicture}
   \end{center}
   \caption{An integral vector in the alternating cone}
   \label{fig1_1}
   \end{figure}

 The \emph{cone of closed alternating walks}, or simply the
 \emph{alternating cone}, ${\mathcal A}(G,{\mathcal C})$ of a
 2-colored graph $(G,\mathcal{C})$ (denoted simply by ${\mathcal
 A}(G)$ when the coloring ${\mathcal C}$ is understood) is defined
 to be the set of all vectors $x = (x(e) : e \in E )$ in ${\mathbb
 R}^E$ satisfying the following system of homogeneous linear
 inequalities:
 \beq \label{3}
 \sum_{e \in E_R(v)} x(e) - \sum_{e \in E_B(v)} x(e) & = & 0, \qquad v \in V, \\
 \label{4}                                           x(e) & \geq &
 0,\qquad e \in E.
 \eeq
 We refer to (\ref{3}) as the \emph{balance condition} at vertex
 $v$. Figure~\ref{fig1_1} illustrates a 2-colored graph together
 with an integral vector in its alternating cone.
 %

 If $G=(V,E)$ is a simple graph, we think of the elements of $E$ as
 2-element subsets of $V$. In this case the 2-colored simple graph
 \emph{associated} to $G$  is the complete graph ${\widehat G} = \left(V,
 {\binom{V}{2}}\right)$, where $e=\{i,j\} \in{\binom{V}{2}}$ is colored
 red if $e \in E$ and colored blue if $e\not \in E$.

 Let $G=(V,E)$ be a graph. A \emph{walk} in $G$ is a sequence
 \beq
 \label{w}
 W &=& (v_0,e_1,v_1,e_2,v_2, \ldots , e_m,v_m),\qquad m \geq 0,
 \eeq
 where $v_i \in V$ for all $i$, $e_j \in E$ for all $j$, and $e_j$
 has endpoints $v_{j-1}$ and $v_j$ for all $j$. We say that $W$ is a
 $v_0$-$v_m$ walk of length $m$. We call $e_1$ the \emph{first} edge
 of $W$ and $e_m$ the \emph{last} edge of $W$. We say that
 $v_1,v_2,\ldots ,v_{m-1}$ are the \emph{internal vertices} of the
 walk $W$. Note that since we are allowing repetitions, the vertices
 $v_0, v_m$ could also be internal vertices. The walk $W^R$ is the
 $v_m$-$v_0$ walk obtained by reversing the sequence (\ref{w}). The
 characteristic vector of the walk $W$ is defined to be $\chi(W) =
 \sum_{i=1}^m \chi(e_i)$.

 The walk $W$ is said to be
 \begin{Lentry}
 \item[closed] when $v_0 = v_m$;

 \item[a trail] when the edges $e_1,\ldots ,e_m$ are distinct;

 \item[a path] when the edges $e_1,\ldots ,e_m$ are distinct and the
 vertices $v_0,\ldots ,v_m$ are distinct;

 \item[a cycle] when $W$ is closed, the edges $e_1,\ldots ,e_m$ are
 distinct, and the vertices $v_0,\ldots ,v_{m-1}$ are distinct.
 \end{Lentry}
 We have defined paths and cycles as special classes of walks.
 However, sometimes it is more convenient to  think of paths and
 cycles as subgraphs, as is done usually. This will be clear from
 the context. If $W_1$ is a $u$-$v$ walk and $W_2$ is a $v$-$w$
 walk, then the \emph{concatenation} of $W_1$ and $W_2$, denoted
 $W_1*W_2$, is the $u$-$w$ walk obtained by walking from $u$
 to $v$ along $W_1$ and continuing by walking from $v$ to $w$
 along $W_2$. Note that if $W_1$ and $W_2$ are trails, then
 $W_1*W_2$ is a trail whenever $W_1$ and $W_2$ have no edges in
 common.

 Now let $(G,\mathcal{C})$ be a 2-colored graph. The walk $W$ in
 (\ref{w}) is said to be
 \begin{Lentry}

 \item[internally alternating] when $\mathcal{C}(e_j) \neq \mathcal{C}(e_{j+1})$ for each
 $j = 1,\ldots, m-1$;

 \item[alternating] when $W$ is internally alternating
 and if $W$ is closed we also have $\mathcal{C}(e_m) \neq
 \mathcal{C}(e_1)$ (note that a walk can be closed and internally
 alternating without being alternating, but if $v_0 \neq v_m$, there
 is no distinction between internally alternating and alternating
 walks and we use the word alternating in this case); a closed
 alternating walk (respectively, trail) is abbreviated as \emph{CAW}
 (respectively, \emph{CAT});

 \item[an even alternating cycle] when $W$ is a cycle of even length and
 $W$ is alternating (Figure~\ref{fig1_2} depicts even alternating cycles and
 their characteristic vectors); an even alternating cycle will also
 be called simply an \emph{alternating cycle};

 \item[an odd internally alternating cycle with base $v_0$] when $W$ is
 a $v_0$-$v_0$ cycle of odd length and $W$ is internally
 alternating (Figure~\ref{fig1_3} depicts odd internally alternating cycles);

 \item[an alternating bicycle] when $W$ is alternating and is of the
 form $W = W_1 * P * W_2 * P^R$, where $W_1, W_2$ are odd internally
 alternating cycles, $P$ is a path between the bases of $W_1$ and
 $W_2$, and the internal vertices of $W_1$, $P$, and $W_2$ are
 disjoint (note that $W_1$ and $W_2$ may have the same base, in
 which case $P$ is empty; Figure~\ref{fig1_4} depicts alternating
 bicycles and their characteristic vectors); clearly $\chi(W) =
 \chi(W_1) + 2\chi(P) + \chi(W_2)$.
 \end{Lentry}
 \begin{figure}[ht]
  \begin{center}
    \begin{pspicture}(-6,-2)(6,2)
      \psset{radius=0.125cm}
      \psset{linewidth=1.8pt}
      \Cnode(-2,0){m0}
      \Cnode(0,0){m1}
      \psset{linewidth=2.5pt}
      \nccurve[angleA=60,angleB=120,linecolor=firebrick1]{m0}{m1}
      \Aput{\textcolor{firebrick1}{$1$}}
      \nccurve[angleA=300,angleB=240,linecolor=skyblue]{m0}{m1}
      \Bput{\textcolor{skyblue}{$1$}}
      \SpecialCoor
      \hspace{4cm}
      \degrees[6]
      \psset{linewidth=1.8pt}
      \Cnode(2;1){n0}
      \Cnode(2;2){n1}
      \Cnode(2;3){n2}
      \Cnode(2;4){n3}
      \Cnode(2;5){n4}
      \Cnode(2;6){n5}
      \psset{linecolor=firebrick1}
      \psset{linewidth=2.5pt}
      \ncline{n0}{n1}\Bput{\textcolor{firebrick1}{$1$}}
      \ncline{n2}{n3}\Bput{\textcolor{firebrick1}{$1$}}
      \ncline{n4}{n5}\Bput{\textcolor{firebrick1}{$1$}}
      \psset{linecolor=skyblue}
      \ncline{n1}{n2}\Bput{\textcolor{skyblue}{$1$}}
      \ncline{n3}{n4}\Bput{\textcolor{skyblue}{$1$}}
      \ncline{n5}{n0}\Bput{\textcolor{skyblue}{$1$}}
    \end{pspicture}
  \end{center}
  \caption{Even alternating cycles}
  \label{fig1_2}
 \end{figure}
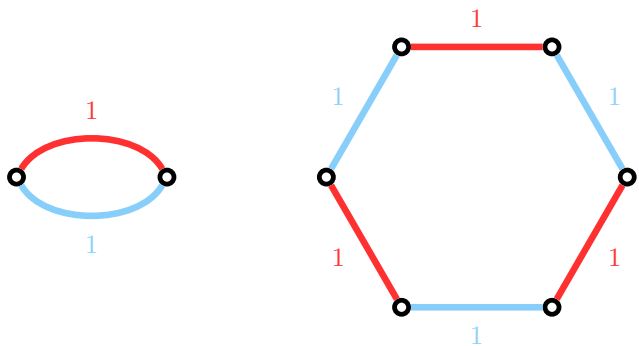
 %
  \begin{figure}[ht]
   \begin{center}
    \begin{pspicture}(2,-2)(18,2.2)
    \psset{radius=0.125cm}
    \psset{linewidth=1.8pt}
    \hspace{6.5cm}
    \SpecialCoor
    \degrees[360]
    \Cnode(1;90){n1}\rput(1.4;90){$v_0$}
    \Cnode(1;210){n2}
    \Cnode(1;330){n3}
    \psset{linecolor=skyblue}
    \psset{linewidth=2.5pt}
    \ncline{n1}{n2}
    \ncline{n3}{n1}
    \psset{linecolor=firebrick1}
    \ncline{n2}{n3}
    \hspace{4.875cm}
    \SpecialCoor
    \degrees[360]
    \psset{linecolor=black}
    \psset{linewidth=1.8pt}
    \Cnode(2;90){m1}\rput(2.4;90){$v_0$}
    \Cnode(2;!360 7 div 90 add){m2}
    \Cnode(2;!720 7 div 90 add){m3}
    \Cnode(2;!3 360 mul 7 div 90 add){m4}
    \Cnode(2;!4 360 mul 7 div 90 add){m5}
    \Cnode(2;!5 360 mul 7 div 90 add){m6}
    \Cnode(2;!6 360 mul 7 div 90 add){m7}
    \psset{linecolor=firebrick1}
    \psset{linewidth=2.5pt}
    \ncline{m1}{m2}
    \ncline{m3}{m4}
    \ncline{m5}{m6}
    \ncline{m7}{m1}
    \psset{linecolor=skyblue}
    \ncline{m2}{m3}
    \ncline{m4}{m5}
    \ncline{m6}{m7}
    \end{pspicture}
  \caption{Odd alternating cycles with base $v_0$}
  \label{fig1_3}
  \end{center}
 \end{figure}
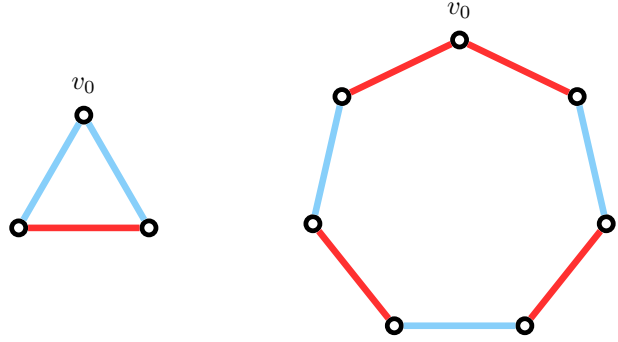
 %
   \begin{figure}[ht]
    \begin{center}
     \begin{pspicture}(-7,-1)(7,1)
      \psset{radius=0.125cm}
      \hspace{-6cm}
      \SpecialCoor
      \degrees[360]
      \psset{linewidth=1.8pt}
      \Cnode(1;0){n1}
      \Cnode(1;72){n2}
      \Cnode(1;144){n3}
      \Cnode(1;216){n4}
      \Cnode(1;288){n5}
      \Cnode(2.1756,0){n6}
      \hspace{3.88cm}
      \Cnode(.7;60){n7}
      \Cnode(.7;180){n8}
      \Cnode(.7;300){n9}
      \psset{linecolor=firebrick1}
      \psset{linewidth=2.5pt}
      \ncline{n1}{n2}\Bput{\textcolor{firebrick1}{$1$}}
      \ncline{n3}{n4}\Bput{\textcolor{firebrick1}{$1$}}
      \ncline{n5}{n1}\Bput{\textcolor{firebrick1}{$1$}}
      \ncline{n6}{n8}\Bput{\textcolor{firebrick1}{$2$}}
      \ncline{n7}{n9}\Aput{\textcolor{firebrick1}{$1$}}
      \psset{linecolor=skyblue}
      \ncline{n2}{n3}\Bput{\textcolor{skyblue}{$1$}}
      \ncline{n4}{n5}\Bput{\textcolor{skyblue}{$1$}}
      \ncline{n1}{n6}\Bput{\textcolor{skyblue}{$2$}}
      \ncline{n7}{n8}\Bput{\textcolor{skyblue}{$1$}}
      \ncline{n8}{n9}\Bput{\textcolor{skyblue}{$1$}}
      \hspace{3cm}
      \psset{linecolor=black}
      \psset{linewidth=1.8pt}
      \Cnode(1;0){m1}
      \Cnode(1;72){m2}
      \Cnode(1;144){m3}
      \Cnode(1;216){m4}
      \Cnode(1;288){m5}
      \Cnode(2.1756,0){m6}
      \Cnode(3.35,0){m7}
      \hspace{5.06cm}
      \Cnode(.7;60){m8}
      \Cnode(.7;180){m9}
      \Cnode(.7;300){m10}
      \psset{linecolor=firebrick1}
      \psset{linewidth=2.5pt}
      \ncline{m1}{m2}\Bput{\textcolor{firebrick1}{$1$}}
      \ncline{m3}{m4}\Bput{\textcolor{firebrick1}{$1$}}
      \ncline{m5}{m1}\Bput{\textcolor{firebrick1}{$1$}}
      \ncline{m6}{m7}\Bput{\textcolor{firebrick1}{$2$}}
      \ncline{m8}{m9}\Bput{\textcolor{firebrick1}{$1$}}
      \ncline{m9}{m10}\Bput{\textcolor{firebrick1}{$1$}}
      \psset{linecolor=skyblue}
      \ncline{m2}{m3}\Bput{\textcolor{skyblue}{$1$}}
      \ncline{m4}{m5}\Bput{\textcolor{skyblue}{$1$}}
      \ncline{m1}{m6}\Bput{\textcolor{skyblue}{$2$}}
      \ncline{m7}{m9}\Bput{\textcolor{skyblue}{$2$}}
      \ncline{m10}{m8}\Bput{\textcolor{skyblue}{$1$}}
     \end{pspicture}
   \caption{Alternating bicycles}
   \label{fig1_4}
  \end{center}
 \end{figure}
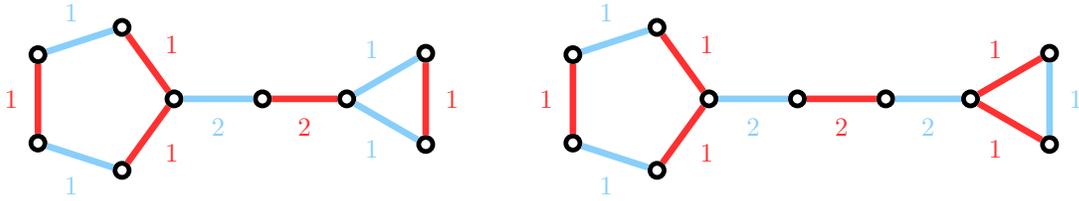

 A CAW $W$ is said to be \emph{irreducible} if $\chi(W)$ cannot be
 written as $\chi(W_1) + \chi(W_2)$ for any CAW's $W_1$ and $W_2$.
 For instance, alternating cycles and bicycles are easily seen to be
 irreducible. Similarly, a CAT $T$ is said to be \emph{irreducible}
 if $\chi(T)$ cannot be written as $\chi(T_1) + \chi(T_2)$ for any
 CAT's $T_1$ and $T_2$. Figure~\ref{fig1_5} depicts an irreducible
 CAW (with the direction of walk indicated by an arrow) and
 Figure~\ref{fig1_6} depicts an irreducible CAT. Irreducibility is
 easily seen.

    \begin{figure}[ht]
     \begin{center}
      \begin{pspicture}(-5,-3.5)(10,1)
       \psset{radius=0.125cm}
       \psset{unit=1.8cm}

       \hspace{-5cm}
       \psset{arrowsize=3pt 2.5}
       \psset{arrowinset=0}
                  \newcommand{\mynclinearrow}
                  {%
                  \nbput[npos=0.1]{\pnode{tail}} \nbput[npos=0.9]{\pnode{head}}%
                  \ncline[linewidth=1pt,linecolor=black,linestyle=dashed]{->}{tail}{head}%
                  }
       \psset{linewidth=1.8pt}
       \Cnode(-.86603,.5){n1}
       \Cnode(-.86603,-.5){n2}
       \Cnode(0,0){n3}
       \Cnode(1,0){n4}
       \Cnode(2,0){n5}
       \Cnode(1.5,-.86603){n6}
       \Cnode(1,-1.73205){n7}
       \Cnode(2,-1.73205){n8}
       \Cnode(2.86603,.5){n9}
       \Cnode(2.86603,-.5){n10}
       \Cnode(3.73205,0){n11}
       \Cnode(4.73205,0){n12}
       \Cnode(5.5981,.5){n13}
       \Cnode(5.5981,-.5){n14}
       \psset{linecolor=firebrick1}
       \psset{linewidth=2.5pt}
       \ncline{n3}{n1}\mynclinearrow
       \ncline{n2}{n3}\mynclinearrow
       \ncline{n4}{n6}\mynclinearrow
       \ncline{n7}{n8}\mynclinearrow
       \ncline{n6}{n5}\mynclinearrow
       \ncline{n10}{n11}\mynclinearrow
       \ncline{n12}{n14}\mynclinearrow
       \ncline{n13}{n12}\mynclinearrow
       \ncline{n11}{n9}\mynclinearrow
       \ncline{n5}{n4}\mynclinearrow
       \psset{linecolor=skyblue}
       \ncline{n1}{n2}\mynclinearrow
       \ncline{n3}{n4}\mynclinearrow
       \ncline{n6}{n7}\mynclinearrow
       \ncline{n8}{n6}\mynclinearrow
       \ncline{n5}{n10}\mynclinearrow
       \ncline{n11}{n12}\mynclinearrow
       \ncline{n14}{n13}\mynclinearrow
       \ncline{n12}{n11}\mynclinearrow
       \ncline{n9}{n5}\mynclinearrow
       \ncline{n4}{n3}\mynclinearrow
      \end{pspicture}
    \caption{An irreducible CAW}
    \label{fig1_5}
    \end{center}
    \end{figure}
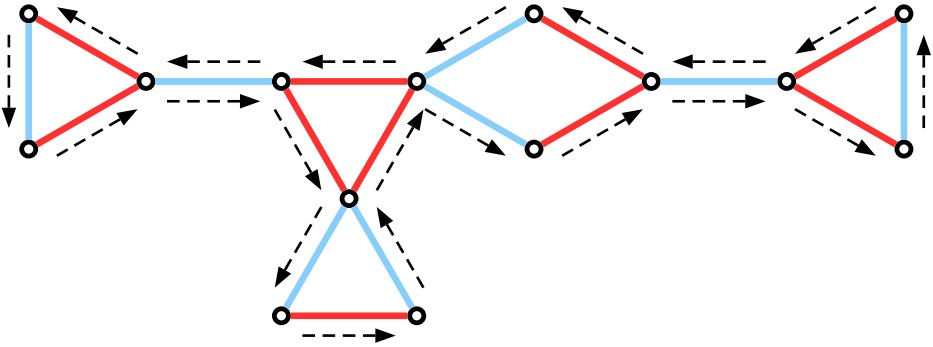
 %
 \begin{figure}[ht]
  \begin{center}
   \begin{pspicture}(-4,-3)(6,3)
    \psset{radius=0.125cm}
    \psset{unit=1.8cm}
    \hspace{-3.5cm}
    \psset{arrowsize=3pt 2.5}
    \psset{arrowinset=0}
               \newcommand{\mynclinearrow}
               {%
               \nbput[npos=0.1]{\pnode{tail}} \nbput[npos=0.9]{\pnode{head}}%
               \ncline[linewidth=1pt,linecolor=black,linestyle=dashed]{->}{tail}{head}%
               }
    \psset{linewidth=1.8pt}
    \Cnode(0,.75){n1}
    \Cnode(-.5,1.61603){n2}
    \Cnode(-1,.75){n3}
    \Cnode(0,-.75){n4}
    \Cnode(-1,-.75){n5}
    \Cnode(-.5,-1.61603){n6}
    \Cnode(1.29904,0){n7}
    \Cnode(1.99,.95107){n8}
    \Cnode(1.99,-.95107){n9}
    \Cnode(3.10806,.5878 ){n10}
    \Cnode(3.10806,-.5878){n11}
    \Cnode(3.60806,1.4538){n12}
    \Cnode(3.60806,-1.4538){n13}
    \Cnode(4.10806,.5878  ){n14}
    \Cnode(4.10806,-.5878 ){n15}
    \psset{linecolor=firebrick1}
    \psset{linewidth=2.5pt}
    \ncline{n2}{n3}\mynclinearrow
    \ncline{n1}{n4}\mynclinearrow
    \ncline{n5}{n6}\mynclinearrow
    \ncline{n4}{n7}\mynclinearrow
    \ncline{n7}{n1}\mynclinearrow
    \ncline{n10}{n8}\mynclinearrow
    \ncline{n9}{n11}\mynclinearrow
    \ncline{n14}{n12}\mynclinearrow
    \ncline{n13}{n15}\mynclinearrow
    \ncline{n11}{n10}\mynclinearrow
    \psset{linecolor=skyblue}
    \ncline{n1}{n2}\mynclinearrow
    \ncline{n3}{n1}\mynclinearrow
    \ncline{n4}{n5}\mynclinearrow
    \ncline{n6}{n4}\mynclinearrow
    \ncline{n8}{n7}\mynclinearrow
    \ncline{n7}{n9}\mynclinearrow
    \ncline{n12}{n10}\mynclinearrow
    \ncline{n10}{n14}\mynclinearrow
    \ncline{n11}{n13}\mynclinearrow
    \ncline{n15}{n11}\mynclinearrow
   \end{pspicture}
   \caption{An irreducible CAT}
   \label{fig1_6}
  \end{center}
 \end{figure}

 Section~\ref{sec2} considers the integral vectors, extreme rays,
 and dimension of the alternating cone. We use a simple alternating
 walk argument to show that the extreme rays of the alternating cone
 are the characteristic vectors of alternating cycles and bicycles,
 the integral vectors in the alternating cone are sums of
 characteristic vectors of irreducible CAW's, and the
 $\{0,1\}$-vectors in the alternating cone are sums of
 characteristic vectors of irreducible CAT's. Using the
 characterization of the extreme rays, we obtain that a simple graph
 $G$ is a threshold graph if and only if $\dim {\mathcal
 A}({\widehat G}) = 0$ (this fact was our original motivation for
 defining the alternating cone). It is well-known that for a simple
 graph  $G$, the property $\dim {\mathcal A}({\widehat G}) = 0$
 (i.e., $G$ being threshold) depends only on the degree sequence of
 $G$. More generally, for any 2-colored graph $(G,\mathcal{C})$, we
 determine $\dim \mathcal{A}(G,\mathcal{C})$ in terms of the red
 degree sequence of $(G,\mathcal{C})$. We then relate this dimension
 to the concept of majorization (following \cite{ap}). Consider the
 set $D(n)$ of all ordered degree sequences $d=(d_1,d_2,\ldots
 ,d_n)$ of simple graphs on $n$ vertices, where $d_1 \geq d_2 \geq
 \cdots \geq d_n$. Partially order $D(n)$ by majorization (the
 definition is recalled in  Section~\ref{sec2} before
 Lemma~\ref{mhl}). It is well-known (see~\cite{mp} and~\cite{rg})
 that the set of maximal elements of this poset is precisely the set
 of ordered degree sequences of threshold graphs. Define a map $A:
 D(n) \rightarrow {\mathbb N}$ by $A(d) = \dim {\mathcal A}(\widehat
 G)$, where $G$ is any simple graph with ordered degree sequence
 $d$. We show that $A$ is an order-reversing map ($d_1 \succeq d_2$
 implies $A(d_1) \leq A(d_2)$). Thus, we can think of $A(d)$ as a
 kind of measure of how non-threshold the degree sequence $d$ is.

 Section~\ref{sec3} is motivated by the following undirected analog
 of Hoffman's circulation problem for directed graphs: let
 $G=(V,E)$, ${\mathcal C} : E \rar \{R,B\}$ be a 2-colored graph.
 Assume that we are given nonnegative integral lower and upper
 bounds $l,u : E \rar {\mathbb N}$ satisfying $l(e) \leq u(e)$ for
 all $e \in E$. We are interested in knowing whether there is a
 vector $y \in {\mathcal A}(G, {\mathcal C}) \cap {\mathbb N}^E$
 satisfying $l(e) \leq y(e)  \leq u(e)$ for all $e \in E$. We use
 residual graph techniques to reduce this problem to the
 \emph{alternating reachability problem}: given distinct vertices
 $s,t$ in a 2-colored graph, is there an alternating $s$-$t$ trail?
 Recall that in the directed case, the circulation problem is
 reduced to the \emph{directed reachability problem}: given distinct
 vertices $s,t$ in a directed graph, is there a directed $s$-$t$
 path? This is solved by a breadth-first search algorithm, which
 either finds a directed $s$-$t$ path or produces an $s$-$t$ cut
 set. In~\cite{BPS2} we give a polynomial-time algorithm to the
 alternating reachability problem generalizing the blossom forest
 algorithm of Edmonds for searching for an augmenting path with
 respect to a matching in a non-bipartite graph. The algorithm
 either finds an alternating $s$-$t$ trail or produces an $s$-$t$ Tutte set
 (which is an obstruction to the existence of an alternating $s$-$t$
 trail. For the definition of a Tutte set, see~\cite{BPS2}).

 Circulations in directed graphs can be thought of in terms of flows.
 For example, the characteristic vector of a directed circuit
 corresponds to a unit of flow along the circuit. Such an
 interpretation is not available in the case of vectors in the
 alternating cone; the irreducible CAW of Figure~\ref{fig1_5} does
 not correspond to a flow  in an intuitive sense. On the other hand,
 the characteristic vector of an irreducible CAT \emph{can} be
 thought of as a unit of flow around the trail. For a 2-colored
 graph $G=(V,E),\; {\mathcal C}: E\rar \{R,B\}$, it is thus natural
 to consider the convex polyhedral cone ${\mathcal T}(G,{\mathcal
 C}) \subseteq {\mathbb R}^E$ generated by  the characteristic
 vectors of the CAT's in $(G,\mathcal{C})$. We call ${\mathcal
 T}(G,{\mathcal C})$ the \emph{cone of closed alternating trails},
 or simply the \emph{trail cone}, of $(G,\mathcal{C})$.

 Consider a CAT in a 2-colored graph. Its characteristic vector
 satisfies the balance condition at every vertex. If we ignore the
 colors, the edge-set of the CAT is a disjoint union of the
 edge-sets of some cycles in the underlying graph. This shows that a
 nonnegative integral combination (that is to say, a linear
 combination with nonnegative integral \emph{coefficients}) of
 characteristic vectors of CAT's satisfies the balance condition at
 every vertex and can be written as a nonnegative integral
 combination of characteristic vectors of cycles in the underlying
 graph $G$. Let ${\mathcal Z}(G)$ denote the cone in ${\mathbb R}^E$
 generated by the characteristic vectors of the cycles in $G$. The
 linear inequalities defining ${\mathcal Z}(G)$ were determined by
 Seymour~\cite{s}. The observation above shows that ${\mathcal
 T}(G,{\mathcal C}) \subseteq {\mathcal A}(G,{\mathcal C}) \cap
 {\mathcal Z}(G)$. In~\cite{BPS2} we prove that ${\mathcal
 T}(G,{\mathcal C}) = {\mathcal A}(G,{\mathcal C}) \cap {\mathcal
 Z}(G)$. The proof uses our solution to the alternating reachability
 problem.

 We remark that in this paper we focus on graph-theoretical aspects
 of the alternating cone and not on algorithmic efficiency. We do
 consider algorithms, but always with a view to obtaining
 graph-theoretical results.

 \section{Extreme Rays and Dimension of the Alternating Cone}\label{sec2}

 A simple graph  $G=(V,E)$ is said to be \emph{threshold} if there
 are real vertex weights $c(v),\;v \in V$ such that every pair
 $e=\{u,v\} \in {\binom{V}{2}}$ satisfies $c(u) + c(v) > 0$ if $e
 \in E$ and $c(u) + c(v) < 0$ if $e \notin E$. Our initial
 motivation for defining the alternating cone was the following
 observation.
 \bt \label{tg} A simple graph $G=(V,E)$ is threshold if and only if
 $\dim {\mathcal A}({\widehat G}) = 0$. \et
 \pf
 Given $e \in {\binom{V}{2}}$, let $\tau (e) = (\tau (e)(v) : v \in V)
 \in {\mathbb R}^V$ denote the incidence vector of $e$, where $\tau
 (e)(v)$ is $1$ if $v$ is an endpoint of $e$ and $0$ otherwise. Let
 ${\mathcal C}_R(G)$ denote the cone in ${\mathbb R}^V$ generated
 by the incidence vectors of the edges $E$, and let ${\mathcal
 C}_B(G)$ denote the cone generated by the incidence vectors of the
 nonedges ${\binom{V}{2}} - E$. If we write (\ref{3}) in matrix
 notation, the columns correspond to the incidence vectors of
 edges  and the negatives of the incidence vectors of nonedges. It
 follows that ${\mathcal A}({\widehat G}) = \{0\}$ if and only if
 ${\mathcal C}_R(G)\, \cap\, {\mathcal C}_B(G) = \{0\}$.

 \noi \textbf{Only if:} Assume that the weights $c(v),\;v \in V$ satisfy the
 defining property of a threshold graph. This means that ${\mathcal
 C}_R(G)$ and ${\mathcal C}_B(G)$ are on opposite sides of the
 hyperplane $\sum_{v \in V} c(v) x(v) = 0$. Hence ${\mathcal C}_R(G)
 \,\cap\, {\mathcal C}_B(G) = \{0\}$.

 \noi \textbf{If:} Suppose ${\mathcal C}_R(G) \,\cap\, {\mathcal
 C}_B(G) = \{0\}$. Then by the separation theorem of convex
 polyhedral cones, there is a hyperplane $\sum_{v \in V} c(v) x(v) =
 0$ such that all nonzero vectors $(p(v): v \in V)  \in {\mathcal
 C}_R(G)$ satisfy $\sum_{v \in V} c(v) p(v) > 0$, and all nonzero
 vectors $(q(v): v \in V)  \in {\mathcal C}_B(G)$ satisfy $\sum_{v
 \in V} c(v) q(v) < 0$. Thus $\{u,v\} \in E$ implies $c(u)+c(v)>0$,
 and $\{u,v\} \notin E$ implies $c(u)+c(v)<0$.
 \epf

 We now determine the extreme rays of the alternating cone.

 \bt
 \label{er}
 Let $G=(V,E),\; {\mathcal C}:E \rar \{R,B\}$ be a 2-colored graph.
 Then
 \be
 \item[\emph{(i)}] the extreme rays of the alternating cone ${\mathcal
 A}(G,{\mathcal C})$ are the characteristic vectors of the
 alternating cycles and bicycles in $(G,\mathcal{C})$;
 \item[\emph{(ii)}] every integral vector in the alternating cone is a nonnegative
 integral combination of the characteristic vectors of irreducible
 CAW's;
 \item[\emph{(iii)}] every $\{0,1\}$-vector in the alternating cone is a
 nonnegative integral combination of the characteristic vectors of
 irreducible CAT's
 \item[\emph{(iv)}] the characteristic vector of an irreducible CAW is
 $\{0,1,2\}$-valued.
 \ee
 \et
 \pf
 (i) Clearly, the characteristic vectors of alternating cycles and
 bicycles are extreme. To show the converse, we will express any
 rational  vector in the alternating cone as a nonnegative rational
 combination of the characteristic vectors of alternating cycles and
 bicycles.

 Let $a=(a(e):e \in E)$ be a nonzero rational  vector in ${\mathcal
 A}(G,{\mathcal C})$. Pick $e_1  \in E$ with $a(e_1)\neq 0$. Without
 loss of generality we may assume that $e_1$ is colored red. Let
 $v_0$ and $v_1$ be the endpoints of $e_1$. Build an alternating
 trail as follows: choose a blue edge $e_2$ incident at $v_1$ with
 $a(e_2)\neq 0$ (this is possible by the balance condition). Let the
 other endpoint of $e_2$ be $v_2$. Now choose a red edge $e_3$
 incident at $v_2$ with $a(e_3) \neq 0$, and so on. At some stage we
 will revisit an already visited vertex. Suppose this happens for
 the first time when we choose edge $e_{k+1}$, i.e., we have built
 an alternating trail
 \beq
 (v_0,e_1,v_1,e_2,v_2, \ldots , e_k,v_k),\qquad k \geq 1,
 \eeq
 where $v_0,v_1,\ldots ,v_k$ are distinct, $\mathcal{C}(e_k) \neq
 \mathcal{C}(e_{k+1})$, and $e_{k+1}$ has endpoints $v_k$ and a
 vertex $u_0 \in\{v_0,v_1,\ldots ,v_{k-1}\}$. Then we have
 found either an alternating cycle $D$, or an odd internally
 alternating cycle $C$ with base $u_0$. In the first case, subtracting an
 appropriate  multiple of $\chi(D)$ from $a$, we obtain another vector
 in the alternating cone whose support is strictly contained in the
 support of $a$. Thus, by induction on the size of the support, we
 are done.

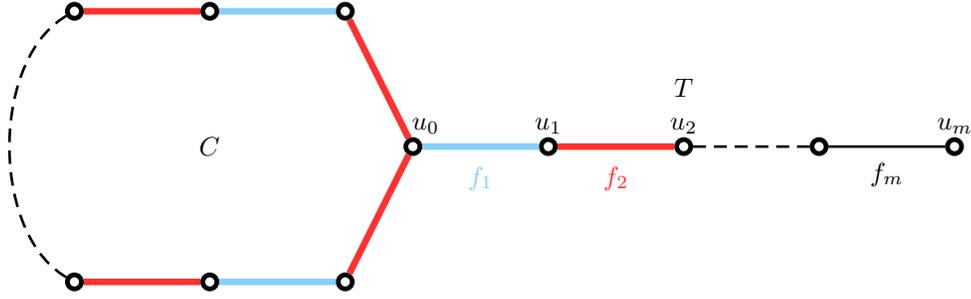
\begin{figure}[ht]
  \begin{center}
   \begin{pspicture}(-6,-2)(7.5,2)
    \hspace{-0.2cm}
    \psset{radius=0.125cm}
    \psset{unit=1.8cm}
    \psset{linewidth=1.8pt}
    \Cnode(-2.5,1){n1}
    \Cnode(-1.5,1){n2}
    \Cnode(-.5,1){n3}
    \Cnode(0,0){n4}
    \Cnode(-.5,-1){n5}
    \Cnode(-1.5,-1){n6}
    \Cnode(-2.5,-1){n7}
    \Cnode(1,0){n8}
    \Cnode(2,0){n9}
    \Cnode(3,0){n10}
    \Cnode(4,0){n11}
    \psset{linecolor=firebrick1}
    \psset{linewidth=2.5pt}
    \ncline{n1}{n2}
    \ncline{n3}{n4}
    \ncline{n4}{n5}
    \ncline{n6}{n7}
    \ncline{n8}{n9}\Bput{\textcolor{firebrick1}{$f_2$}}
    \psset{linecolor=skyblue}
    \ncline{n5}{n6}
    \ncline{n2}{n3}
    \ncline{n4}{n8}\Bput{\textcolor{skyblue}{$f_1$}}
    \psset{linecolor=black}
    \ncline[linewidth=1 pt,linestyle=dashed]{n9}{n10}
    \ncline[linewidth=1 pt]{n10}{n11}\Bput{\textcolor{black}{$f_m$}}
    \ncarc[linewidth=1 pt,arcangle=60,linestyle=dashed]{n7}{n1}
    \rput (-1.5,0){$C$}
    \uput{5pt}[65](0,0){$u_0$}
    \uput{5pt}[90](1,0){$u_1$}
    \uput{5pt}[90](2,0){$u_2$}
    \uput{5pt}[90](4,0){$u_m$}
    \uput{1pt}[90](2,.35){$T$}
   \end{pspicture}
  \end{center}
  \caption{Illustrating the proof of Theorem
  \protect\ref{er} (i)}
  \label{fig2_1}
 \end{figure}

 In the second case, extend $C$ to an alternating trail $C*T$ as
 follows: $T$ starts with $T = (u_0,f_1,\ldots)$, where $f_1$ is an
 edge incident with $u_0$ and satisfies $a(f_1)\neq 0$ and
 $\mathcal{C}(f_1) \neq \mathcal{C}(e_{k+1})$. Let the other endpoint
 of $f_1$ be $u_1$. Now add to $T$ an edge $f_2$ incident with $u_1$
 and satisfying $a(f_2)\neq 0$ and $\mathcal{C}(f_2) \neq
 \mathcal{C}(f_1)$, and so on. At some stage we will revisit an
 already visited vertex of the trail $C*T$. Suppose this happens for
 the first time when we choose edge $f_{m+1}$, i.e., we have
 \beq
 T & =& (u_0,f_1,u_1,f_2,u_2, \ldots , f_m,u_m),\qquad m \geq 1,
 \eeq
 where $u_0,u_1,\ldots ,u_m$ are distinct, none of $\{u_1,\ldots
 ,u_m\}$ is on $C$ (see Figure~\ref{fig2_1}), $\mathcal{C}(f_{m+1})
 \neq \mathcal{C}(f_m)$, one endpoint of $f_{m+1}$ is $u_m$, and the
 other endpoint $v$ of $f_{m+1}$ is either in $T$ or is a vertex of
 $C$ different from $u_0$. Two cases arise:
 \\
 \noi \textbf{Case (a):} $v$  is a vertex of $T$ (see
 Figure~\ref{fig2_2}). We have found either an alternating cycle or
 an alternating bicycle, and we are done by induction on the size of
 the support, as in the previous paragraph.
 \\
 \noi \textbf{Case (b):} $v$ is not a vertex of $T$ (see
 Figure~\ref{fig2_3}). In this case, the edges $f_1,\ldots f_{m+1}$
 together with an appropriate portion of $C$ determine an
 alternating cycle, and we are done.
 \\

 \begin{figure}[ht]
  \begin{center}
   \hspace{-1cm}
   \begin{pspicture}(-6,-2)(9,2)
    \psset{radius=0.125cm}
    \psset{unit=1.8cm}
    \psset{linewidth=1.8pt}
    \Cnode(-2.5,1){n1}
    \Cnode(-1.5,1){n2}
    \Cnode(-.5,1){n3}
    \Cnode(0,0){n4}
    \Cnode(-.5,-1){n5}
    \Cnode(-1.5,-1){n6}
    \Cnode(-2.5,-1){n7}
    \Cnode(1,0){n8}
    \Cnode(2,0){n9}
    \Cnode(3,0){n10}
    \Cnode(4,0){n11}
    \Cnode(5,0){n12}
    \psset{linecolor=firebrick1}
    \psset{linewidth=2.5pt}
    \ncline{n1}{n2}
    \ncline{n3}{n4}
    \ncline{n4}{n5}
    \ncline{n6}{n7}
    \ncline{n8}{n9}\Bput{\textcolor{firebrick1}{$f_2$}}
    \psset{linecolor=skyblue}
    \ncline{n5}{n6}
    \ncline{n2}{n3}
    \ncline{n4}{n8}\Bput{\textcolor{skyblue}{$f_1$}}
    \psset{linecolor=black}
    \ncline[linewidth=1 pt,linestyle=dashed]{n9}{n10}
    \ncline[linewidth=1 pt,linestyle=dashed]{n10}{n11}
    \ncline[linewidth=1 pt]{n11}{n12}\Bput{\textcolor{black}{$f_m$}}
    \ncarc[linewidth=1 pt,arcangle=90]{n12}{n10}
    \ncarc[linewidth=1 pt,arcangle=60,linestyle=dashed]{n7}{n1}
    \rput (-1.5,0){$C$}
    \uput{5pt}[65](0,0){$u_0$}
    \uput{5pt}[90](1,0){$u_1$}
    \uput{5pt}[90](2,0){$u_2$}
    \uput{5pt}[90](3,.025){$v$}
    \uput{5pt}[90](5,0){$u_m$}
    \uput{1pt}[90](2.5,.35){$T$}
    \uput{15pt}[270](4,-0.415){$f_{m+1}$}
   \end{pspicture}
  \end{center}
  \caption{Illustrating the proof of Theorem \protect\ref{er} (i) Case (a)}
  \label{fig2_2}
 \end{figure}
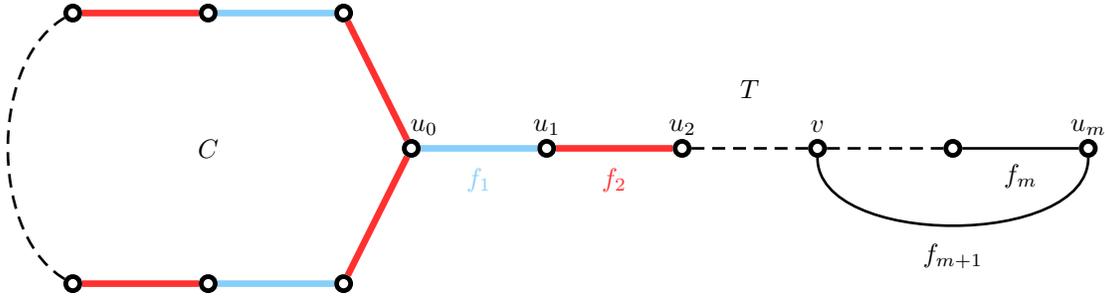

 \begin{figure}[ht]
  \begin{center}
   \begin{pspicture}(-6,-3.6)(8,2)
    \psset{radius=0.125cm}
    \psset{unit=1.8cm}
    \psset{linewidth=1.8pt}
    \Cnode(-2.5,1){n1}
    \Cnode(-1.5,1){n2}
    \Cnode(-.5,1){n3}
    \Cnode(0,0){n4}
    \Cnode(-.5,-1){n5}
    \Cnode(-1.5,-1){n6}
    \Cnode(-2.5,-1){n7}
    \Cnode(-3.08,0){nn1}
    \Cnode(1,0){n8}
    \Cnode(2,0){n9}
    \Cnode(3,0){n10}
    \Cnode(4,0){n11}
    \psset{linecolor=firebrick1}
    \psset{linewidth=2.5pt}
    \ncline{n1}{n2}
    \ncline{n3}{n4}
    \ncline{n4}{n5}
    \ncline{n6}{n7}
    \ncline{n8}{n9}\Bput{\textcolor{firebrick1}{$f_2$}}
    \psset{linecolor=skyblue}
    \ncline{n5}{n6}
    \ncline{n2}{n3}
    \ncline{n4}{n8}\Bput{\textcolor{skyblue}{$f_1$}}
    \psset{linecolor=black}
    \ncline[linewidth=1 pt,linestyle=dashed]{n9}{n10}
    \ncline[linewidth=1 pt]{n10}{n11}\Bput{\textcolor{black}{$f_m$}}
    \ncarc[linewidth=1 pt,arcangle=29,linestyle=dashed]{n7}{nn1}
    \ncarc[linewidth=1 pt,arcangle=29,linestyle=dashed]{nn1}{n1}
    \ncarc[linewidth=1 pt,arcangle=110,linestyle=solid]{n11}{nn1}
    \rput (-1.5,0){$C$}
    \uput{5pt}[65](0,0){$u_0$}
    \uput{5pt}[90](1,0){$u_1$}
    \uput{5pt}[90](2,0){$u_2$}
    \uput{5pt}[90](4,0){$u_m$}
    \uput{1pt}[90](2,.35){$T$}
    \uput{1pt}[270](1,-1.9){$f_{m+1}$}
   \end{pspicture}
  \end{center}
  \caption{Illustrating the proof of Theorem \protect\ref{er} (i) Case (b)}
  \label{fig2_3}
 \end{figure}
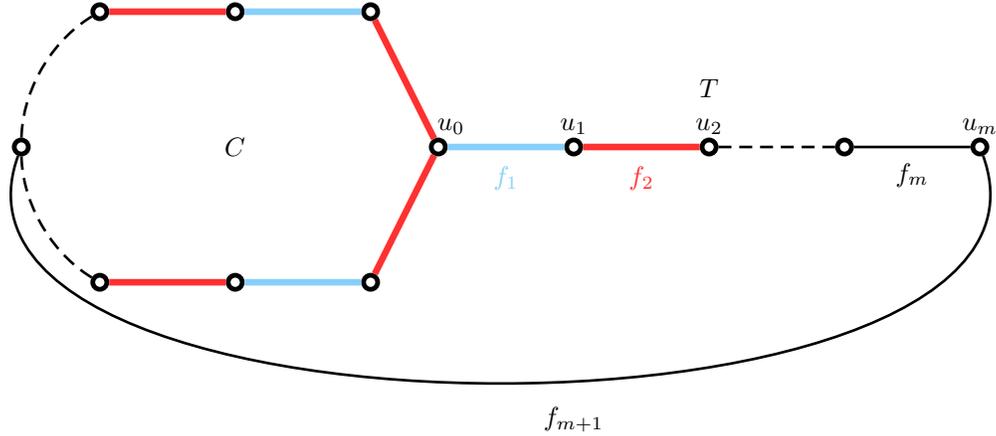

 Essentially the same argument as given above appears in~\cite{hip}
 (in the context of edges and non-edges).

 \noi (ii) Let $a=(a(e):e \in E)$ be a nonzero  integral vector in
 the alternating cone. Pick an edge $e_1$ with $a(e_1)\neq 0$ and
 with end points $v_0$ and $v_1$. Assume that we have an internally
 alternating walk
 \beq W &=& (v_0,e_1,v_1,e_2,v_2, \ldots , e_m,v_m),\qquad m \geq 1,
 \eeq
 with $\chi(W) \leq a$ (we can always start with the walk
 $(v_0,e_1,v_1)$). We show below that either we can extend $W$, or
 else there is a CAW (and hence an irreducible CAW) through $e_1$.
 Since we cannot extend indefinitely because of the condition
 $\chi(W) \leq a$, we are done. The following cases arise.
 \\
 \noi \textbf{Case (a):} $v_m\neq v_0$. Then $\chi(W)$ does not
 satisfy the balance condition at $v_m$, but $a$ does, and since
 $\chi(W) \leq a$ and $a$ is integral, we can find an edge $e_{m+1}$
 incident at $v_m$ with $\mathcal{C}(e_{m+1}) \neq \mathcal{C}(e_m)$
 such that $\chi(W)(e_{m+1}) < a(e_{m+1})$. Extend $W$ by adding
 $e_{m+1}$ and the other end point of $e_{m+1}$.
 \\
 \noi \textbf{Case (b):} $v_m = v_0$ and $\mathcal{C}(e_1) =
 \mathcal{C}(e_m)$. We can extend $W$ just as in case~(a).
 \\
 \noi \textbf{Case (c):} $v_m = v_0$ and $\mathcal{C}(e_1) \neq
 \mathcal{C}(e_m)$. In this case $W$ is a CAW.

 \noi (iii) This is a special case of (ii): if $a$ is a
 $\{0,1\}$-vector, then the CAW's in (ii) must be CAT's.

 \noi (iv) Consider a walk $W$ as in~(\ref{w}). This assigns a
 direction of traversal to each edge; for instance, the edge $e_2$
 is traversed from $v_1$ to $v_2$. The direction of traversal may be
 different for two occurrences of the same edge. However, if a CAW
 $W$ traverses an edge three or more times, then two of these
 directions must be the same, and this can be used to write
 $W=W_1*W_2$ for two positive length CAW's $W_1$ and $W_2$, so $W$
 is not irreducible. \epf

 As a corollary of Theorem~\ref{er}, we derive the following well-known
 characterization of threshold graphs.
 \bco A simple graph $G$ is not threshold if and only if $\widehat
 G$ contains an alternating cycle of length 4.
 \eco
 \pf
 \textbf{If:} Suppose $\{ \{i,j\},\{j,k\},\{k,l\},\{l,i\} \}$ is an
 alternating 4-cycle in $\widehat G$ with the pairs
 $\{i,j\},\{k,l\}$ red and the other two pairs blue. Assume that $G$
 is threshold with vertex weights $c(v),\;v \in V$ satisfying the
 defining property. Since $\{i,j\},\{k,l\}$ are red, we have
 $c(i)+c(j) >0$, $c(k)+c(l)>0$ and therefore $c(i)+c(j)+c(k)+c(l) >
 0$. Similarly, since $\{j,k\},\{l,i\}$ are blue, we have
 $c(i)+c(j)+c(k)+c(l) < 0$, a contradiction.

 \noi \textbf{Only if:} Since $G$ is not threshold, by Theorem~\ref{tg}
 ${\mathcal A}(\widehat G)$ has an extreme ray, which is an alternating cycle
 or an alternating bicycle by Theorem~\ref{er}. Suppose that this
 extreme ray is an alternating cycle of length greater than $4$.
 There is a chord of $\widehat G$ that splits this cycle into two
 even cycles (since $\widehat G$ is complete). Regardless of the
 color of the chord, one of these two cycles is alternating.
 Repeating this argument, we obtain an alternating cycle of length
 4.

 Now consider an extreme ray that is an alternating bicycle
 $W=W_1*P*W_2*P^R$. Let $u$ and $v$ be the bases of $W_1$ and $W_2$.
 Let $u'$ (respectively, $v'$) be any vertex of $W_1$ (respectively,
 $W_2$) different from $u$ (respectively, $v$). Consider the edge
 $\{u',v'\}$ of $\widehat G$. $W_1$ determines two alternating
 $u'$-$u$ paths, and one of them starts with an edge having color
 different from that of $\{u',v'\}$. Call this alternating path
 $P_1$. Similarly, using $W_2$, choose an alternating $v'$-$v$ path
 $P_2$ that starts with an edge having color different from that of
 $\{u',v'\}$. We now  have the alternating cycle
 $P_1*P*P_2^R*(v',\{v',u'\},u')$, and we can use the argument of the
 preceding paragraph.
 \epf

 We now give a formula for the dimension of the alternating cone of
 a 2-colored graph.

 A connected graph is said to be \emph{odd unicyclic} if it contains
 precisely one cycle, and that cycle has odd length. In other words,
 an odd unicyclic graph is obtained from a tree by adding a new edge
 between two nonadjacent vertices of the tree so that the cycle
 created has odd length. A graph is a \emph{pseudo forest} if each
 component of the graph is either acyclic or odd unicyclic. Pseudo forests
 are to be distinguished from $1$-forests, which are graphs whose connected
 components have at most one cycle, even or odd.  The motivation for
 studying $1$-forests is combinatorial while pseudo forests have a linear
 algebraic origin (see Theorem~\ref{qft} below).
 Recall from the proof of Theorem~\ref{tg} that for an
 edge $e$, the vector $\tau(e) \in {\mathbb R}^V$ is the incidence
 vector of $e$, which is 1 in the two coordinates indexed by the
 endpoints of $e$, and is 0 elsewhere. The \emph{incidence matrix}
 of a graph is the matrix whose columns are the incidence vectors of
 the edges. For a proof of the following result see~\cite{gks}.
 \bt
 \label{qft}
 For a graph $G=(V,E)$ and a set $X \subseteq E$, the set $\{ \tau(e) :
 e \in X \}$ is linearly independent in ${\mathbb R}^V$ if and only
 if the graph $(V,X)$ is a pseudo forest. In particular, the rank of
 the incidence matrix of $G$ is equal to $\#V - \mbox{ number of
 bipartite components of } G$.
 \et
 For a graph $G=(V,E)$ and an integer sequence $d=(d(v) :
 v \in V)$, we use the notation
 \[{\mathcal K}(d) = \{{\mathcal C}:E\rightarrow \{R,B\} : \mbox{the red degree of }
 v \mbox{ is } d(v) \mbox{ for all } v \in V\},\]
 i.e., ${\mathcal K}(d)$ denotes the set of all 2-colorings of $G$
 having red degree sequence $d$.
 \bl
 \label{cl}
 Consider the 2-colored graph $G=(V,E)$ with a coloring ${\mathcal
 C} \in {\mathcal K}(d)$, and let $e \in E$. If $(G,\mathcal{C})$
 has a CAW through $e$, then for each ${\mathcal C'} \in {\mathcal
 K}(d)$, $(G,\mathcal{C'})$ has a CAW through $e$.
 \el
 \pf
 Let ${\mathcal C'} \in {\mathcal K}(d)$ and consider the spanning subgraph
 $G'=(V,E')$ of $G$, where $E'$ consists of all the edges where
 $\mathcal{C}$ and $\mathcal{C'}$ disagree.
 %
 %

 Since the red degrees (and thus also the blue degrees) in $G$ agree
 under ${\mathcal C}$ and ${\mathcal C'}$, it follows that for each
 $v \in V$, the red degree of $v$ in $(G',\mathcal{C'})$ is equal to
 the blue degree of $v$ in $(G',\mathcal{C'})$. Thus the all-1
 vector in ${\mathbb R}^{E'}$ is balanced in $(G',\mathcal{C'})$,
 i.e., is in $\mathcal{A}(G',\mathcal{C'})$. It follows from
 Theorem~\ref{er}(iii) that for each $e \in E'$, $(G',\mathcal{C'})$
 has a CAT through $e$, and therefore so does $(G,\mathcal{C'})$.

 Now let $W$ be a CAW through $e$ in $(G,\mathcal{C})$. If
 ${\mathcal C}(e) \neq {\mathcal C'}(e)$, then we already know that
 $(G,\mathcal{C'})$ has a CAT through $e$, and we are done. So we
 may assume that ${\mathcal C}(e) = {\mathcal C'}(e)$. We will
 transform $W$ into a CAW $W'$ through $e$ in $(G,\mathcal{C'})$.
 Let $f$ be an edge in $W$ with endpoints $u$ and $v$. If ${\mathcal
 C}(f)={\mathcal C'}(f)$, we do nothing. If ${\mathcal C}(f) \neq
 {\mathcal C'}(f)$, then $f \in E'$ and $(G,\mathcal{C'})$ has a CAT
 through $f$. Dropping $f$ from this CAT, we obtain  $u$-$v$
 alternating trail $P$ in $(G,{\mathcal C'})$ whose first and last
 edges have the color ${\mathcal C}(f)$. We drop $f$ from $W$ and
 substitute the trail $P$ in its place. Doing this for every edge
 $f$ in $W$ with ${\mathcal C}(f)\neq {\mathcal C'}(f)$, we obtain a
 CAW $W'$ through $e$ in $(G,\mathcal{C'})$. \epf

 For a graph $G$ and an integer sequence $d$, let $E_d$ be
 the set of all edges $e$ of $G$ such that some 2-coloring in
 $\mathcal{K}(d)$ has a CAW through $e$ (equivalently by
 Lemma~\ref{cl}, all 2-colorings in $\mathcal{K}(d)$ have a CAW
 through $e$).

 \bt
 \label{dt}
 Let $G=(V,E)$ be a graph and ${\mathcal C}$ a 2-coloring of $G$
 with red degree sequence $d$. Then
 \[\dim {\mathcal A}(G,{\mathcal C}) = \#{E_d} - \#V +
 b(V,E_d),\]
 where $b(V,E_d)$ denotes the number of bipartite components of the graph $(V,E_d)$.
 \et
 \pf
 An edge $e \in E$ is said to be \emph{inessential} if $x(e)=0$ for
 all $x \in {\mathcal A}(G,{\mathcal C})$. From Theorem~\ref{er}(ii)
 it follows that $e$ is inessential if and only if $(G,\mathcal{C})$
 has no CAW through $e$. From basic polyhedral theory it now follows
 that $\dim {\mathcal A}(G,{\mathcal C})$ is equal to the nullity
 (i.e., number of columns minus rank) of the $\#V\times \#E_d$
 vertex-edge incidence matrix of the graph $(V,E_d)$. The expression
 for the dimension now follows from Theorem~\ref{qft}.
 \epf

 From Theorem~\ref{dt}, $\dim {\mathcal A}(G,{\mathcal C})$ depends
 only on $G$ and the red degree sequence of $\mathcal{C}$. In the case of
 the associated 2-colored graphs of simple graphs we can say more.
 \bl \label{sdac} Let $G_1$ and $G_2$ be simple graphs with degree
 sequences $d_1$ and $d_2$. If $d_1$ is a rearrangement of $d_2$ (so
 in particular $G_1$ and $G_2$ have the same number of vertices),
 then $\dim {\mathcal A}(\widehat G_1) = \dim {\mathcal A}(\widehat
 G_2)$. \el
 \pf
 Suppose that the permutation $\pi :V \rar V$ rearranges $d_1$ into
 $d_2$. The result follows from the fact that $\pi$ is an automorphism
 of the complete graph $\left(V,\binom{V}{2}\right)$.
 \epf

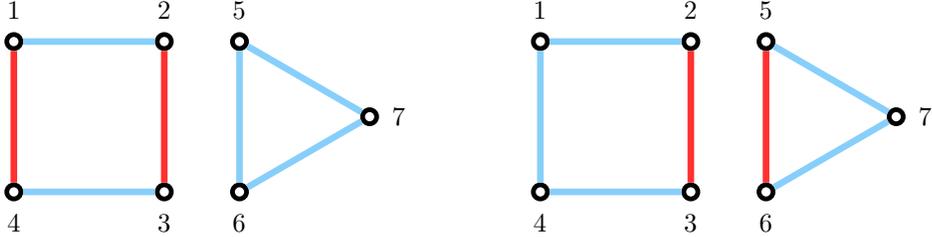
\begin{figure}[ht]
  \begin{center}
   \begin{pspicture}(-11,-2)(1,2)
    \hspace{-21.7cm}
    \psset{radius=0.125cm}
    \psset{linewidth=1.8pt}
    \Cnode(0,-1){n1}
    \Cnode(2,-1){n2}
    \Cnode(2,1){n3}
    \Cnode(0,1){n4}
    \Cnode(3,-1){n5}
    \Cnode(3,1){n6}
    \Cnode(4.73,0){n7}
    \uput{5pt}[90](0,1.125){$1$}
    \uput{5pt}[90](2,1.125){$2$}
    \uput{5pt}[270](2,-1.125){$3$}
    \uput{5pt}[270](0,-1.125){$4$}
    \uput{5pt}[270](3,-1.125){$6$}
    \uput{5pt}[90](3,1.125){$5$}
    \uput{5pt}[0](4.855,0){$7$}
    \psset{linecolor=firebrick1}
    \psset{linewidth=2.5pt}
    \ncline{n1}{n4}
    \ncline{n3}{n2}
    \psset{linecolor=skyblue}
    \ncline{n2}{n1}
    \ncline{n3}{n4}
    \ncline{n5}{n6}
    \ncline{n7}{n6}
    \ncline{n7}{n5}
    \psset{linecolor=black}
    \psset{linewidth=1.8pt}
    \Cnode(7,-1){mn1}
    \Cnode(9,-1){mn2}
    \Cnode(9,1){mn3}
    \Cnode(7,1){mn4}
    \Cnode(10,-1){mn5}
    \Cnode(10,1){mn6}
    \Cnode(11.73,0){mn7}
    \psset{linecolor=firebrick1}
    \psset{linewidth=2.5pt}
    \ncline{mn2}{mn3}
    \ncline{mn5}{mn6}
    \psset{linecolor=skyblue}
    \ncline{mn2}{mn1}
    \ncline{mn3}{mn4}
    \ncline{mn1}{mn4}
    \ncline{mn7}{mn6}
    \ncline{mn7}{mn5}
    \uput{5pt}[90](7,1.125){$1$}
    \uput{5pt}[90](9,1.125){$2$}
    \uput{5pt}[270](9,-1.125){$3$}
    \uput{5pt}[270](7,-1.125){$4$}
    \uput{5pt}[270](10,-1.125){$6$}
    \uput{5pt}[90](10,1.125){$5$}
    \uput{5pt}[0](11.855,0){$7$}
   \end{pspicture}
  \end{center}
  \caption{2-colored graphs with the same red degree sequence:
  1-dimensional alternating cone (left) and 0-dimensional
  alternating cone (right)}
  \label{fig2_4}
 \end{figure}

 Lemma~\ref{sdac} fails for 2-colored graphs that are not complete:
 \bex
 Figures~\ref{fig2_4} depicts two 2-colorings of a graph on the
 vertex set $\{1,2,\ldots ,7\}$ whose red degree sequences are
 permutations of each other (via the permutation $\pi$ that fixes
 $2,3,7$ and exchanges $1$ with $5$ and $4$ with $6$). However, it
 is easily seen that the dimensions of the alternating cones of the
 2-colored graphs are 1 and 0, respectively. The permutation $\pi$
 is not an automorphism of the underlying graph.
 \eex

 %


 We now relate the dimension of the alternating cone to the concept
 of majorization. We begin with a few definitions.

 Let $a=(a(1),\ldots ,a(n))$ and $b=(b(1),\ldots ,b(n))$ be real
 sequences of length $n$. Denote the $i$-th largest component of $a$
 (respectively, $b$) by $a[i]$ (respectively, $b[i]$). We say that
 $a$ \emph{majorizes} $b$, denoted by $a \succeq b$, if
 \[\sum_{i=1}^k a[i]  \geq \sum_{i=1}^k b[i],\qquad k=1,\ldots ,n,\]
 with equality for $k=n$. The majorization is \emph{strict}, denoted
 by $a \succ b$, if at least one of the inequalities is strict, namely
 if $a$ is not a permutation of $b$. We recall a fundamental
 lemma about majorization in integer sequences, called Muirhead's
 lemma. If $a=(a(1),\ldots ,a(n))$ is a sequence and there exist $i$ and $j$
 such that  $a(i) \geq a(j) + 2$,
 then the following operation is called a \emph{unit
 transformation from $i$ to $j$ on $a$}: subtract 1 from $a(i)$ and
 add 1 to $a(j)$. Clearly, if $b$ is obtained from $a$ by a sequence
 of unit transformations, then $a\succ b$. The converse is also true
 for integer sequences.

 \bt[Muirhead Lemma] \label{mhl} If $a$ and $b$ are integer
 sequences and $a\succ b$, then some permutation of $b$ can be
 obtained from $a$ by a sequence of unit transformations. \et
 For a proof see \cite{mp,mo}.

 \bt
 \label{dml}
 Let $G_1 = (V,E_1)$ and $G_2 = (V,E_2)$ be simple graphs with
 degree sequences $d_1$ and $d_2$. If $d_1 \succeq d_2$, then $\dim
 {\mathcal A}(\widehat G_1) \leq \dim {\mathcal A}(\widehat G_2)$.
 \et
 \pf
 If $d_2$ is a rearrangement of $d_1$, the result follows from
 Lemma~\ref{sdac}, so we may assume that $d_1 \succ d_2$. By
 Muirhead's lemma some permutation $d_2'$ of $d_2$ can be obtained
 from $d_1$ by a finite sequence $d_1 \succ d \succ \cdots \succ
 d_2'$ of unit transformations. We will show that $d$ is the degree
 sequence of a simple graph $G$ satisfying $\dim {\mathcal
 A}(\widehat G_1)  \leq \dim {\mathcal A}(\widehat G)$. By
 Lemma~\ref{sdac} and induction on the number of unit
 transformations, the result will follow.

\begin{figure}[ht]
  \begin{center}
   \begin{pspicture}(-4,-3)(4,2)
    \hspace{-2.95cm}
    \psset{radius=0.125cm}
    \hspace{-3cm}
    \Cnode(0,-1){n1}
    \Cnode(2,-1){n2}
    \Cnode(2,1){n3}
    \Cnode(0,1){n4}
    \Cnode(4,-1){m1}
    \Cnode(6,-1){m2}
    \Cnode(6,1){m3}
    \Cnode(4,1){m4}
    \ncline{n3}{n4}
    \ncline{n4}{n2}
    \ncline[linestyle=dashed]{n1}{n3}
    \ncline[linestyle=dashed]{n1}{n2}
    \ncline{m3}{m4}
    \ncline{m1}{m2}
    \ncline[linestyle=dashed]{m1}{m3}
    \ncline[linestyle=dashed]{m2}{m4}
    \uput{5pt}[90](0,1.125){$i$}
    \uput{5pt}[90](2,1.125){$l$}
    \uput{5pt}[270](0,-1.125){$j$}
    \uput{5pt}[270](2,-1.125){$k$}
    \uput{5pt}[90](4,1.125){$i$}
    \uput{5pt}[90](6,1.125){$l$}
    \uput{5pt}[270](4,-1.125){$j$}
    \uput{5pt}[270](6,-1.125){$k$}
    \rput(1,-2.5){$G_1$}
    \rput(5,-2.5){$G$}
   \end{pspicture}
  \end{center}
  \caption{Illustrating the proof of Theorem \protect\ref{dml}}
  \label{fig2_6}
 \end{figure}
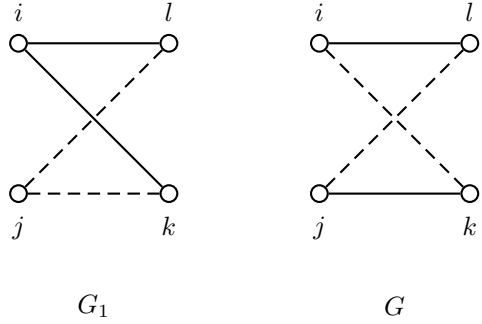

 For notational convenience, let $V = \{1,2,\ldots,n\}$, and
 suppose $d$ is obtained from $d_1$ by a unit transformation from
 $i$ to $j$, so that $d_1(i) \geq d_1(j) + 2$. This implies that there
 exist distinct vertices $k,l \neq i,j$ such that $\{i,k\},
 \{i,l\}$ are edges of $G_1$ and $\{j,k\}, \{j,l\}$ are not. Let $G$
 be the graph with degree sequence $d$ obtained from $G_1$ by
 dropping the edge $\{i,k\}$ and adding the edge $\{j,k\}$ (see
 Figure~\ref{fig2_6}).

 Consider the 2-colored graphs $\widehat{G}_1$ and $\widehat{G}$
 with red degree sequences $d_1$ and $d$, respectively. We now show
 that $E_{d_1}  \subseteq E_d$. By Theorem~\ref{dt}, it will then
 follow that $\dim {\mathcal A}(\widehat G_1)  \leq \dim {\mathcal
 A}(\widehat G)$, since $\#(E_d - E_{d_1})  \geq b(V,E_{d_1}) -
 b(V,E_d)$. This last inequality can be seen as follows: start with
 the graph $(V,E_{d_1})$ and add the edges $e \in E_d - E_{d_1}$ one
 at a time. If $e$ connects two components $C_1$ and $C_2$, the
 number of bipartite components decreases by one or stays the same,
 according as $C_1$ and $C_2$ are both bipartite or not; if $e$
 connects two vertices in the same component $C$, the number of
 bipartite components stays the same if $C$ is nonbipartite, and it
 stays the same or decreases by one if $C$ is bipartite according to
 the parity of (any of the) cycles created by $e$.

 Suppose $\{u,v\} \in E_{d_1}$. If $\{u,v\}$ is one of the pairs
 $\{i,k\},\{j,k\}$ that changed status by going from $G_1$ to $G$,
 then Figure~\ref{fig2_6} depicts an alternating 4-cycle through
 $\{u,v\}$ in $\widehat G$, and thus $\{u,v\} \in E_d$ and we are
 done. So we may assume that $\{u,v\}$ is not one of these two
 pairs.

 Since $\{u,v\} \in E_{d_1}$, $\widehat G_1$ has a CAW $W$ through
 $\{u,v\}$. Replace every occurrence of
 \begin{gather*}
 \ldots,i,\{i,k\},k,\ldots,
 \\
 \ldots,k,\{i,k\},i,\ldots,
 \\
 \ldots,k,\{k,j\},j,\ldots,
 \\
 \ldots,j,\{j,k\},k,\ldots,
 \end{gather*}
 in $W$ by (respectively)
 \begin{gather*}
 \ldots,i,\{i,l\},l,\{l,j\},j,\{j,k\},k,\ldots,
 \\
 \ldots,k,\{k,j\},j,\{j,l\},l,\{l,i\},i,\ldots,
 \\
 \ldots,k,\{k,i\},i,\{i,l\},l,\{l,j\},j,\ldots,
 \\
 \ldots,j,\{j,l\},l,\{l,i\},i,\{i,k\},k,\ldots,
 \end{gather*}
 keeping all other edges in $W$ fixed. This yields a CAW
 through $\{u,v\}$ in $G$, and thus $\{u,v\}  \in E_d$.
 \epf

 As stated in the introduction, Theorem~\ref{dml} defines an
 order-reversing map $A: D(n) \rar {\mathbb N}$, which maps the
 degree sequence of a simple graph $G$ to the dimension of the
 alternating cone of the associated 2-colored graph $\widehat{G}$.
 Given $d=(d(1),\ldots ,d(n))  \in D(n)$, there is a well-known
 algorithm working only with the numbers $d(1),\ldots ,d(n)$ to
 determine whether $A(d)=0$ (see \cite{mp}). Motivated by this, we
 ask whether there is an algorithm working only with the numbers
 $d(1),\ldots ,d(n)$ for computing $A(d)$.

 \section{Intersection of the Alternating Cone with a Box}\label{sec3}

 Assume that we are given a 2-colored graph $G=(V,E)$, and for each
 $e \in E$ nonnegative integers $l(e),u(e)$ with $l(e) \leq u(e)$.
 We ask if there is a rational vector $x \in {\mathcal
 A}(G,{\mathcal C})$ with $l(e) \leq x(e) \leq u(e)$ for all $e \in
 E$. The next theorem restricts the search to half-integral $x$.
 \bt
 \label{hi}
 Let $G=(V,E)$, ${\mathcal C}: E\rar \{R,B\}$ be a 2-colored graph,
 and $l,u : E\rightarrow {\mathbb N}$ maps with $l(e) \leq u(e)$ for
 all $e \in E$. If there exists a rational vector $x \in {\mathcal
 A}(G,{\mathcal C})$ with $l(e) \leq x(e) \leq u(e)$ for all $e \in
 E$, then there exists an integral $y \in {\mathcal A}(G,{\mathcal
 C})$ with $2l(e) \leq y(e) \leq 2u(e)$ for all $e \in E$.
 \et
 \pf
 We use elementary  polyhedral theory. Since by assumption a feasible
 solution exists, there exists a basic feasible solution $\overline{x}
 \in {\mathcal A}(G,{\mathcal C})$, with $l \leq \overline{x} \leq u$. In
 our case a basic feasible solution is obtained as follows. First
 choose a pseudo forest $(V,X)$ such that the columns corresponding
 to $X$ form a basis of the column space of the vertex-edge
 incidence matrix of $G$. For each $e \in E - X$ we have $\overline{x}(e)
 = l(e)$ or $\overline{x}(e) = u(e)$. Now solve for the remaining
 $\overline{x}(e)$, $e \in X$ using the balance condition at every node.
 Since $l,u$ are integral and the determinant of the incidence
 matrix of an odd cycle is $\pm2$, the half-integrality of $\overline{x}$
 easily follows, and $y = 2\overline{x}$ is as required.
 \epf

 Motivated by Theorem~\ref{hi}, we want to improve
 half-integrality to integrality, so we are led to the following
 problem. Let a 2-colored graph $G=(V,E)$, ${\mathcal C}: E\rar
 \{R,B\}$ and bounds $l,u : E \rightarrow {\mathbb N}$ be given. For
 $f \in \ac$, an edge $e$ is called \emph{feasible w.r.t.\ $f$} if
 $l(e) \leq f(e) \leq u(e)$, and $f$ itself is called
 \emph{feasible} if every edge is feasible w.r.t.\ $f$,
 \emph{infeasible} otherwise. We ask if there is a feasible vector
 $f \in \ac \cap \mn^E$. We now reduce this problem to the problem of
 finding a CAT through a given edge in a 2-colored graph. This latter
 problem is easily reduced to the alternating reachability problem.

 Let $f \in \ac \cap \mn^E$, $f$ not necessarily feasible. The
 \emph{residual 2-colored graph $G(f)=(V,E(f))$ of $f$ w.r.t.\
 $l,u$} is defined as follows. We take four disjoint copies $E_1,
 E_2, E_3, E_4$ of $E$, and denote the copy of $e \in E$ in $E_i$ by
 $e_i$, $i=1,\ldots,4$. For each $e \in E$, we place $e_1$ in $E(f)$
 with the color $\cc(e)$ when $f(e) \leq u(e) - 1$, place $e_2$ in
 $E(f)$ with the color $\cc(e)$ when $f(e) \leq u(e) - 2$, place
 $e_3$ in $E(f)$ with the color opposite $\cc(e)$ when $f(e) \geq
 l(e) + 1$, and place $e_4$ in $E(f)$ with the color opposite
 $\cc(e)$ when $f(e) \geq l(e) + 2$.

 Suppose that $G(f)$ has a CAT $T$. We extend the characteristic
 vector $\chi(T)$ by adding zero components at all elements of $E_1
 \cup E_2 \cup E_3 \cup E_4 - E(f)$. By \emph{augmenting $f$ along
 $T$} we mean replacing $f$ with $f_T$ given by
 \[f_T(e) = f(e) + \chi(T)(e_1) + \chi(T)(e_2) - \chi(T)(e_3) -
 \chi(T)(e_4), \qquad e \in E.\]
 Note that $f_T  \in \ac \cap \mn^E$, and that in replacing $f$ with
 $f_T$, feasible edges remains feasible, the infeasible edges of $T$
 move ``in the right direction'', i.e., become feasible or move
 closer to feasibility, and of course the edges out of $T$ remain
 unchanged.
 \bt
 \label{at} Suppose $f \in \ac\cap {\mn}^E$ is infeasible, but
 $\ac\cap {\mn}^E$ has a feasible vector. Then for each $e \in E$,
 \be
 \item [\emph{(i)}] if $f(e) < l(e)$, then  $G(f)$ has a CAT through $e_1$;
 \item [\emph{(ii)}] if $f(e) > u(e)$, then  $G(f)$ has a CAT through $e_3$.
 \ee
 \et
 \pf
 We define a 2-colored subgraph $G'(f)=(V,E'(f))$ of $G(f)$ by
 letting $E'(f) = E(f) \cap (E_1 \cup E_3)$ and restricting the
 2-coloring of $G(f)$ to $G'(f)$.

 Let $g \in \ac \cap \mn^E$ be a feasible vector. We define $h : E'(f)\rar \mn$ as
 follows: for $e_1 \in E'(f)$,
 \[ h(e_1) = \left\{ \ba{ll}
                      0 & \mbox{if } g(e) - f(e) < 0, \\
                      g(e) - f(e) & \mbox{if } g(e) - f(e) \geq 0,
                      \ea
              \right. \]
 and for $e_3 \in E'(f)$,
 \[ h(e_3) = \left\{ \ba{ll}
                      0 & \mbox{if } g(e) - f(e) > 0,
                      \\
                      -(g(e) - f(e)) & \mbox{if } g(e) - f(e) \leq 0.
                      \ea
              \right. \]
 It is easy to check that $h$ is an integral vector in
 the alternating cone of $G'(f)$.

 \noi (i) Assume that $e \in E$ with $f(e) < l(e)$. Then $e_1 \in
 E'(f)$ and $h(e_1) > 0$ (since $g$ is feasible). By
 Theorem~\ref{er}(ii), $h$ can be written as a sum of characteristic
 vectors of irreducible CAW's in $G'(f)$, and thus $G'(f)$ has an
 irreducible CAW $W$ through $e_1$. By Theorem~\ref{er}(iv),
 $\chi(W)$ is $\{0,1,2\}$-valued. Suppose $\chi(W)(a_1) = 2$ for
 some $a_1$ (respectively, $\chi(W)(a_3) = 2$ for some $a_3$). Then
 $g(a) - f(a) \geq 2$ (respectively, $f(a) - g(a) \geq 2$). Since
 $g$ is feasible, $a_2 \in E(f)$ (respectively, $a_4 \in E(f)$), and
 consequently $a_1 \in E(f)$ (respectively, $a_3 \in E(f)$) by the
 definition of $E(f)$. We then consider $W$ as a subset of $E(f)$
 and replace the double occurrence of $a_1$ (respectively, $a_3$) in
 $W$ by a single occurrence of $a_1$ and of $a_2$ (respectively, of
 $a_3$ and of $a_4$). Doing this for all repeated edges in $W$
 transforms it into a CAT in $G(f)$ through $e_1$.

 \noi (ii) Similar to (i).
 \epf

 The problem of finding a CAT through a given edge $e$ in an
 edge-colored graph can be reduced to the alternating
 reachability problem as follows: let $e$ have endpoints  $s$ and $t$. Remove $e$
 from the graph, add two new vertices $s'$ and $t'$, add two new edges
 with color $\mathcal{C}(e)$, one  between  $s'$ and $s$ and one
 between $t'$ and $t$. Clearly the new graph has an alternating $s'$-$t'$
 trail if and only if the original graph has a CAT through $e$.

 We can now use the following familiar scheme to look for a feasible
 integral vector. Start with an integral balanced $f:E\rar \mn$, for
 example $f = 0$. If $f$ is infeasible, construct $G(f)$. Pick an
 edge $e \in E$ with $f(e)< l(e)$ (respectively, $f(e) > u(e)$), and
 find a CAT $T$ through $e_1$ (respectively, $e_3$) in $G(f)$ if one
 exists (using the alternating reachability algorithm in~\cite{BPS2}), then
 augment $f$ along $T$. As noted in the proof of Theorem~\ref{at},
 $f_T$ is integral and balanced, feasible edges remain feasible, and
 in addition each infeasible edge in $T$, in particular $e$, has
 either become feasible or has moved closer to feasibility. Replace
 $f$ with $f_T$ and repeat. Since we are working with integral
 vectors, either we terminate with a feasible integral vector in
 time bounded by the total infeasibility, or else at some stage
 $G_f$ has no CAT through $e_1$ (respectively, $e_3$), in which case
 no feasible integral vector exists by Theorem~\ref{at}.

 As stated in the introduction, in this paper we are not dealing
 with efficiency issues but only with the graph-theoretic aspects of
 the alternating cone. Our discussion motivates the alternating
 reachability problem and the problem of determining the linear
 inequalities defining the trail cone. These two problems are
 considered in~\cite{BPS2}.

\noindent
{\bf Acknowledgement:} We thank the referees for their constructive
suggestions that led to an improvement in the exposition.

 \end{document}